\documentclass{proc-l}
\usepackage{amsfonts}

\newtheorem{theorem}{Theorem}
\theoremstyle{plain}

\newtheorem{corollary}{Corollary}

\newtheorem{remark}{Remark}

\numberwithin{equation}{section}

\title{Probabilistic aspects of Al-Salam--Chihara polynomials}
\author{W{\l }odzimierz Bryc}
\address{Department of Mathematics \\
University of Cincinnati\\
P.O. Box 210025\\
Cincinnati, OH 45221--0025}
\thanks{\noindent Research partially supported by NSF
grant \#INT-0332062.}
 \email{ Wlodzimierz.Bryc@UC.edu}
\author{Wojciech Matysiak}
\address{Faculty of Mathematics and Information Science\\
Warsaw University of Technology\\
pl. Politechniki 1\\
 00-661 Warszawa, Poland}
\email{wmatysiak@elka.pw.edu.pl}
\author{Pawe\l\ J. Szab\l owski}
\address{ Faculty of Mathematics and Information Science\\
Warsaw University of Technology\\ pl. Politechniki 1\\ 00-661
Warszawa, Poland}
\email{pszablowski@elka.pw.edu.pl}
\date{April 11, 2003; Revised: November 31, 2003}


\begin{document}

\begin{abstract}
We solve the connection coefficient problem between the Al-Salam--Chihara
polynomials and the $q$-Hermite polynomials, and we use the resulting
identity to answer a question from probability theory. We also derive the
distribution of some Al-Salam--Chihara polynomials, and compute determinants
of related Hankel matrices.
\end{abstract}

\maketitle


\section{Introduction and main identity}

The aim of the paper is to point out the connection of Al-Salam--Chihara
polynomials with a regression problem in probability, and to use it to give
a new simple derivation of their density. Our approach exploits identity
(\ref{MI}) below, which connects the Al-Salam--Chihara polynomials to the
continuous $q$-Hermite polynomials. This
connection is more direct and elementary but less general than  the technique of attachment exploited
in \cite[Section 2]{Berg-Ismail96}.
We also compute determinants of Hankel
matrices with entries which are linear combinations of the $q$-Hermite
polynomials.

The Al-Salam--Chihara polynomials were introduced in
\cite{AlSalam-Chihara76SIAM} and their weight function was found
in \cite{Askey-Ismail84MAMS}. We are interested in the
renormalized Al-Salam--Chihara polynomials $\{p_{n}(x|q,a,b)\}$,
which are defined by the following three term recurrence relation
\begin{equation}
p_{n+1}\left( x\right) =\left( x-aq^{n}\right) p_{n}\left( x\right) -\left(
1-bq^{n-1}\right) [n]_{q}p_{n-1}\left( x\right) \ (n\geq 0),  \label{*}
\end{equation}%
with the usual initial conditions
$p_{-1}(x)=0$, $p_{0}(x)(x) =1$.
Here, we use the standard notation
\begin{eqnarray*}
{[n]_{q}} &=&1+q+\dots +q^{n-1}, \\
{[n]_{q}!} &=&[1]_{q}[2]_{q}\dots \lbrack n]_{q}, \\
\left[
\begin{array}{c}
n \\
k%
\end{array}%
\right] _{q} &=&\frac{[n]_{q}!}{[n-k]_{q}![k]_{q}!},
\end{eqnarray*}%
with the usual conventions $[0]_{q}=0,[0]_{q}!=1$.

For $|q|<1$, their generating function
\begin{equation*}
f(t,x|q,a,b)=\sum_{n=0}^\infty \frac{t^n}{[n]_q!}p_n(x|q,a,b)
\end{equation*}
is given by
\begin{equation}  \label{**}
f(t,x|q,a,b)=\prod_{k=0}^\infty \frac{1-(1-q)at q^k+(1-q)b t^2 q^{2k}}{%
1-(1-q)xtq^k+(1-q)t^2q^{2k}};
\end{equation}
compare \cite[(3.6) and (3.10)]{Askey-Ismail84MAMS}.

The corresponding (renormalized) continuous $q$-Hermite polynomials $%
H_n(x|q)=p_n(x|q,0,0)$ satisfy the three term recurrence relation
\begin{equation}  \label{H}
H_{n+1}(x)=x H_n(x)-[n]_qH_{n-1}(x).
\end{equation}
For $|q|<1$ their generating function $\phi(t,x|q)=\sum_{n=0}^\infty
\frac{t^n}{[n]_q!}H_n(x|q)$ is
\begin{equation}  \label{***}
\phi(t,x|q)= \prod_{k=0}^\infty
\left(1-(1-q)xtq^k+(1-q)t^2q^{2k}\right)^{-1}.
\end{equation}
Of course, these are well known special cases of (\ref{*}) and (\ref{**}),
see \cite[(2.11) and (2.12)]{ISV87} which we state here for further reference.

We will also use polynomials $\{B_{n}(x|q)\}$ defined by the three term
recurrence relation
\begin{equation}
B_{n+1}\left( x\right) =-q^{n}xB_{n}\left( x\right)
+q^{n-1}[n]_{q}B_{n-1}\left( x\right) \ (n\geq 0)  \label{B}
\end{equation}%
with the usual initial conditions $B_{-1}=0,B_{0}=1$.
 These polynomials are related to the $q$-Hermite polynomials by
\begin{equation}
B_{n}(x|q)=\left\{
\begin{array}{ll}
i^{n}q^{n(n-2)/2}H_{n}(i\sqrt{q}\,x|\frac{1}{q}) & \mbox{ if $q>0$} \\
(-1)^{n(n-1)/2}|q|^{n(n-2)/2}H_{n}(-\sqrt{|q|}\,x|\frac{1}{q}) & \mbox{ if $q<0$}%
\end{array}%
\right.,  \label{bnaH}
\end{equation}%
and have been studied in \cite{Askey-89}, \cite{Ismail-Masson}.
 Their generating function
$\psi (t,x|q)=\sum_{n=0}^{\infty }\frac{t^{n}}{[n]_{q}!}B_{n}(x|q)$
is given by
\begin{equation}
\psi (t,x|q)=\prod_{k=0}^{\infty }\left(
1-(1-q)xtq^{k}+(1-q)t^{2}q^{2k}\right).  \label{****}
\end{equation}

We now point out the mutual relationship between the Al-Salam--Chihara polynomials $%
\{p_n(x|q,a,b)\}$ and the polynomials $\left\{ H_{n}\left( x|q\right)
\right\} _{n\geq 0}$ and $\left\{ B_{n}\left( x|q\right) \right\} _{n\geq 0}$.
\begin{theorem}
\label{MainT} For all $a,c,q\in \mathbb{C}$, 
$c\neq 0$, and $n\geq 1$ we have
\begin{equation}
p_{n}(x|q,a,b)=\sum_{k=1}^{n}\left[
\begin{array}{c}
n \\
k%
\end{array}%
\right] _{q}c^{n-k}B_{n-k}(\frac{a}{c}|q)\left( H_{k}(x|q)-c^{k}H_{k}(\frac{a%
}{c}|q)\right)  \label{MI},
\end{equation}
where $b=c^2$.
\end{theorem}

\begin{proof}
From the recurrence relations (\ref{*}), (\ref{H}), and (\ref{B}),
it is clear that $p_{n}(x|q,a,b)$, $H_{n}(x|q)$, and $B_{n}(x|q)$
are given by polynomial expressions in variable $q$.
The $q$-binomial coefficient $\left[
\begin{array}{c}
n \\
k
\end{array}%
\right] _{q}$ is also a polynomial in $q$.
Therefore, we see that identity (\ref{MI}) is equivalent to a polynomial identity in
variable $q\in \mathbb{C}$. Hence it is enough to prove that (\ref{MI})
holds true for all $|q|<1$. When $|q|<1$, inspecting (\ref{**}), (\ref{***}%
), and (\ref{****}) we notice that for $b=c^2$ we have
\begin{equation}
f(t,x|q,a,b)=\psi (tc,a/c|q)\phi (t,x|q),  \label{g1}
\end{equation}%
and
\begin{equation}
\psi (t,x|q)\phi (t,x|q)=1.  \label{g2}
\end{equation}%
Therefore,
\begin{equation*}
f(t,x|q,a,b)=1+\psi (ct,a/c|q)\left( \phi (t,x|q)-\phi (ct,a/c|q)\right)
,
\end{equation*}%
which is valid for all small enough $|t|$. Comparing the coefficients at $%
t^{n}$ for $n\geq 1$ and taking into account that $H_{k}(x|q)-c^{k}H_{k}(%
\frac{a}{c}|q)=0$ for $k=0$, we get (\ref{MI}).
\end{proof}

\begin{remark}
One could split (\ref{MI}) into the following two separate identities, which
are implied by (\ref{g1}) and (\ref{g2}) respectively:
\begin{equation}\label{*1}
\forall n\geq 0: \ p_{n}(x|q,a,c^{2})=\sum_{k=0}^{n}\left[
\begin{array}{c}
n \\
k%
\end{array}%
\right] _{q}c^{n-k}B_{n-k}(\frac{a}{c}|q) H_{k}(x|q),
\end{equation}
\begin{equation}\label{*2}
\forall n\geq 1: \ \sum_{k=0}^{n} 
\left[\begin{array}{c}n\\k\end{array}\right]_qB_{n-k}(x|q)H_{k}\left(
x|q\right) =0.
\end{equation}
Formula (\ref{*1}) is a renormalized inverse of formula
\cite[(4.7)]{IRS99}, who express the $q$-Hermite polynomials as
linear combinations of Al-Salam--Chihara polynomials. Formula
(\ref{*2}) resembles \cite[(2.28)]{Carlitz56}, who considers
$q$-Hermite polynomials of the form $h_n(x|q)=\displaystyle
\sum_{k=0}^n\left[
\begin{array}{c}
n \\ k\end{array}\right]_{q}x^k$, paired with  $b_n(x|q)=h_n(x|1/q)$.
\end{remark}

\section{Probabilistic aspects}

Quadratic regression questions in the paper \cite{Bryc01} lead to the
problem of determining all 
probability distributions $\mu$ which are defined indirectly by the
relationships
\begin{equation}  \label{warunkowy_q_Hermit}
\int H_n(x|q)\mu(dx)=\rho^nH_n(y|q) \ n=1,2,\dots,
\end{equation}
where $y,\rho, q \in\mathbb{R}$ are fixed parameters, and $\left\{
H_{n}\right\} _{n\geq 0}$ is the family of the $q$-Hermite polynomials.

Our next result shows that this problem can be solved using the
Al-Salam--Chihara polynomials.

\begin{theorem}
\label{T2} If $\mu =\mu (dx|\rho ,y)$ satisfies (\ref{warunkowy_q_Hermit}),
then its orthogonal polynomials are Al-Salam--Chihara polynomials $%
\{p_{n}(x|q,a,b)\}$ with $a=\rho y$, $b=\rho ^{2}$.
\end{theorem}

\begin{proof}
Recall that $H_{n}(x|q)=p_{n}(x|q,0,0)$. Thus if $\rho =0$ then (\ref%
{warunkowy_q_Hermit}) implies that $\int p_{n}(x|q,a,b)\mu (dx)=0$ for all $%
n=1,2,\dots $. Suppose now that $\rho \neq 0$. Combining (\ref{MI}) with (%
\ref{warunkowy_q_Hermit}) we get
\begin{equation*}
\int p_{n}(x|q,a,b)\mu (dx)=\sum_{k=0}^{n}\left[
\begin{array}{c}
n \\
k%
\end{array}%
\right] _{q}\rho ^{n-k}B_{n-k}(y|q)\int \left( H_{k}(x|q)-\rho
^{k}H_{k}(y|q)\right) \mu (dx)=0
\end{equation*}%
for all $n=1,2,\dots $. Since $\{p_n\}$ satisfy a three step recurrence,
this implies
$\int p_k(x)p_{n}(x)\mu(dx)=0$ for all $0\leq k<n$.
\end{proof}

Next we answer an unresolved case from \cite{Bryc01}.

\begin{corollary} Fix $q>1, y\in \mathbb{R}$. Let ${\mathcal R}_q=\{1,1/q,1/q^2,\dots,1/q^n,\dots,0\}$.
\begin{itemize}
\item[(i)] If  $\rho^2\not\in {\mathcal R}_q$ then (\ref{warunkowy_q_Hermit}) has no
probabilistic solution $\mu$.
\item[(ii)] If  $\rho^2\in {\mathcal R}_q$ is non-zero,  then the probabilistic solution of
(\ref{warunkowy_q_Hermit}) exists, and is a discrete measure 
supported on $1+\log_q 1/\rho^2$ points.

\end{itemize}
\end{corollary}

\begin{proof}
Suppose that $\mu $ is positive and solves (\ref{warunkowy_q_Hermit}).
Therefore
its monic orthogonal polynomials satisfy the three term recurrence relation 
\begin{equation}
p_{n+1}(x)=(x-\rho yq^{n})p_{n}(x)-(1-\rho ^{2}q^{n-1})p_{n-1}(x).  \label{+}
\end{equation}%
For a positive non-degenerate measure $\mu_y(dx)$, and $n\geq 1$ we have
\begin{equation}\label{+++}
\int p_{n}^2(x)\mu_y(dx)=(1-\rho ^{2}q^{n-1})\int p_{n-1}^2(x)\mu_y(dx).\end{equation}
 if  $\rho^2\not\in {\mathcal R}_q$ then $(1-\rho ^{2}q^{n-1})\ne 0$ for all $n$.
Since $\int p_{0}^2(x)\mu_y(dx)>0$, this shows that
$\int p_{n}^2(x)\mu_y(dx)>0$ for all $n\geq 0$.  But then
 the coefficients $1-\rho ^{2}q^{n-1}$
must be non-negative for all $n$, which is false. This proves (i).

To conclude the proof it remains to notice that if $\rho^2=1/q^m$  then
from
(\ref{+}) and (the proof of) Favard's theorem, see \cite[Theorem II.1.5]{Freud},
 it follows that the solution of (\ref{warunkowy_q_Hermit}) is given by a measure supported on the roots
of the polynomial $p_{m+1}$.
 Indeed, (\ref{+}) implies that the polynomial $p_{m+2}$ is divisible by
$p_{m+1}$. Therefore, $p_{m+1}$ is the common factor of
all polynomials  $\{p_k:k\geq m+1\}$.
It is also known, see \cite[Theorem I.2.2]{Freud},
that $p_{m+1}$ has exactly $m+1$ distinct real
roots   $x_1,\dots,x_{m+1}$. Thus, any measure $\mu(dx)=\sum \lambda_j \delta_{x_j}$
supported on the roots
of polynomial
$p_{m+1}$ satisfies  $\int p_{m+1+k}\mu(dx)=0$. Solving the remaining $m+1$ equations
$\int p_0 \mu(dx)=1$, and
$\int p_k(x)\mu(dx)=0,k=1,2,\dots,m$ for $\lambda_j$, we get a
 measure that solves
(\ref{warunkowy_q_Hermit}). This measure is non-negative as
 the coefficients at the third term in the recurrence (\ref{+})
are non-negative for $n=1,\dots,m$, see \cite[page 58]{Freud}.
\end{proof}

From Theorem \ref{T2} it follows that if the solution of (\ref%
{warunkowy_q_Hermit}) exists, then it is given by the distribution of the
Al-Salam-Chihara polynomials. The distribution of the Al-Salam-Chihara
polynomials is derived in \cite[Chapter 3]{Askey-Ismail84MAMS}. However, in
\cite[Proposition 8.1]{Bryc01} we found the solution of (\ref%
{warunkowy_q_Hermit}) which relies solely on the facts about the $q$-Hermite
polynomials. We repeat the latter argument here, and then use it to
re-derive the distribution of the corresponding Al-Salam--Chihara
polynomials.

\begin{corollary}
\label{C2} If $\rho ,q, y\in \mathbb{R}$ are such that $|\rho |<1$,$|q|<1$, and $y^2(1-q)<4$, then
the probabilistic solution of (\ref{warunkowy_q_Hermit}) is given by the
absolutely continuous measure $\mu $ with the density on $x^{2}<4/(1-q)$
given by
\begin{equation*}
\frac{\sqrt{1-q}}{2\pi \sqrt{4-(1-q)x^{2}}}\prod_{k=0}^{\infty }\frac{%
(1-\rho ^{2}q^{k})\left( 1-q^{k+1}\right) \left(
(1+q^{k})^{2}-(1-q)x^{2}q^{k}\right) }{(1-\rho ^{2}q^{2k})^{2}-(1-q)\rho
q^{k}(1+\rho ^{2}q^{2k})xy+(1-q)\rho ^{2}(x^{2}+y^{2})q^{2k}}
\end{equation*}
\end{corollary}

\begin{proof}
The distribution of the $q$-Hermite polynomials $H_{n}(x|q)$ is supported on
$x^{2}<4/(1-q)$ with the density \begin{equation*}
f_{H}(x)=\frac{\sqrt{1-q}}{2\pi \sqrt{4-(1-q)x^{2}}}\prod_{k=0}^{\infty
}\left( (1+q^{k})^{2}-(1-q)x^{2}q^{k}\right) \prod_{k=0}^{\infty
}(1-q^{k+1}),
\end{equation*}%
see \cite[(2.15)]{ISV87}.
Moreover, since $|H_n(x)|\leq C_q(n+1)(1-q)^{-n/2}$ when $x^2,y^2\leq 4/(1-q)$, the
series
\begin{equation}
g_H(x,y,\rho )=\sum_{n=0}^{\infty }\frac{\rho ^{n}}{[n]_{q}!}H_{n}(x)H_{n}(y)
\label{PM}
\end{equation}
converges uniformly and defines the Poisson-Mehler kernel
which is given by
\begin{equation}\label{q-M}
g_H(x,y,\rho )=\prod_{k=0}^{\infty }\frac{(1-\rho ^{2}q^{k})}{(1-\rho
^{2}q^{2k})^{2}-(1-q)\rho q^{k}(1+\rho ^{2}q^{2k})xy+(1-q)\rho
^{2}(x^{2}+y^{2})q^{2k}};
\end{equation}%
this is the renormalized version of the well known result, see eg.
\cite[(2.2)]{IS88} who consider the
$q$-Hermite polynomials given by
$\{(1-q)^{n/2}H_n(2x/\sqrt{1-q}|q)\}$ instead of our $\{H_n(x|q)\}$.

%
%
%

Since (\ref{H}) implies that $\int H_{n}^{2}(x|q)f_{H}(x)dx=[n]_{q}!$, it
follows from (\ref{PM}) that
\begin{equation*}
\int_{-2/\sqrt{1-q}}^{2/\sqrt{1-q}}H_{n}(x|q)g_H(x,y,\rho )f_{H}(x)\,dx=\rho
^{n}H_{n}(y).
\end{equation*}
\end{proof}

\begin{corollary}\label{C3}
If $q,a,b\in\mathbb{R}$ are such that $|q|<1$, $0< b<1$, and $a^2(1-q)<4b$, then the distribution of the
Al-Salam-Chihara polynomials $\{p_n(x|q,a,b)\}$ is absolutely continuous
with the density on  $x^2<4/(1-q)$ given by
\begin{equation*}
\frac{\sqrt{1-q}}{2\pi\sqrt{4-(1-q)x^2}} \prod_{k=0}^\infty \frac{%
(1-bq^k)\left(1-q^{k+1}\right)\left((1+q^{k})^2-(1-q)x^2q^{k}\right)} {%
(1-bq^{2k})^2-(1-q)a q^k(1+bq^{2k})x+(1-q)(bx^2+a^2)q^{2k}}.
\end{equation*}
\end{corollary}

\begin{proof}
By Theorem \ref{T2}, the distribution of polynomials $p_{n}$ solves (\ref%
{warunkowy_q_Hermit}) with $\rho =\sqrt{b},y=a/\rho $. Thus the formula
follows from Corollary \ref{C2}.
\end{proof}
\begin{remark}
Iterating (\ref{warunkowy_q_Hermit}) we see that the measure
corresponding to the parameter $\rho_1\rho_2$ instead of $\rho$ is given
by
\begin{eqnarray}\label{Chapman}
\mu (\cdot|\rho_1\rho_2 ,x)=\int \mu(\cdot|\rho_1,y) \mu (dy|\rho_2 ,x).
\end{eqnarray}
For $\left\vert q\right\vert <1, |\rho|<1$ the density of $\mu$ is given in
Corollary \ref{C2}, hence after simplifying common factors and
substitution $x=2\xi /\sqrt{1-q},$ $y=2\eta /\sqrt{1-q},$
$z=2\zeta /\sqrt{1-q}$, the relationship (\ref{Chapman}) takes the
following form
\begin{eqnarray*}
&&\int_{-1}^{1}\prod_{k=0}^{\infty }\frac{(1-\rho_1
^{2}q^{k})\left( 1-q^{k+1}\right) \left( (1+q^{k})^{2}-4\eta
^{2}q^{k}\right) }{(1-\rho_1 ^{2}q^{2k})^{2}-4\rho_1 q^{k}(1+\rho_1
^{2}q^{2k})\eta \zeta +4\rho_1 ^{2}(\eta
^{2}+\zeta ^{2})q^{2k}} \\
&&\times \prod_{k=0}^{\infty }\frac{(1-\rho_2
^{2}q^{k})}{(1-\rho_2
^{2}q^{2k})^{2}-4\rho_2 q^{k}(1+\rho_2 ^{2}q^{2k})\xi \eta +4\rho_2
^{2}(\eta ^{2}+\xi ^{2})q^{2k}}\frac{d\eta}{2\pi
\sqrt{1-\eta ^{2}}}  \\
&=&\prod_{k=0}^{\infty }\frac{(1-\rho_1^2\rho_2^2 q^{k})}{(1-\rho_1^2\rho_2^2q^{2k})^{2}-4\rho_1\rho_2q^{k}
(1+\rho_1^2\rho_2^2q^{2k})\xi \zeta +4\rho_1^2\rho_2^2(\zeta ^{2}+\xi ^{2})q^{2k}}.
\end{eqnarray*}
\end{remark}
\section{Determinants of Hankel matrices}

In this section we are interested in calculating the determinants of the
Hankel matrices
\begin{equation*}
M_n=[m_{i+j}]_{i,j=0,\ldots n-1},
\end{equation*}%
where $m_{i} =\int x^{i} \mu(dx)$ are the moments of certain (perhaps
signed) measure $\mu$. It is well known that for positive measures we must
have $\det M_n\geq 0$, and that these determinants can be read out from the
three term recurrence for the corresponding monic orthogonal polynomials.

Consider first
 the moments $m_k(y)=\int x^k \mu(dx)$ of the (perhaps signed)
measure $\mu=\mu_{y,\rho}$, which solves (\ref{warunkowy_q_Hermit}). Then $%
m_k(y)$ are polynomials of degree $k$ in variable $y$, and can be written as
follows. Let $a_{n,2i}$, $i\leq \left\lfloor n/2\right\rfloor $ be the
coefficients in the expansion of monomial $x^n$ into the $q$-Hermite
polynomials,
\begin{equation*}
x^{n}=\sum_{i=0}^{\left\lfloor n/2\right\rfloor }a_{n,2i}H_{n-2i}\left(
x|q\right) ,\ n\geq 0.
\end{equation*}
Then
\begin{equation*}
m_n(y)=\int x^{n}d\mu \left( x\right) =\sum_{k=0}^{\left\lfloor
n/2\right\rfloor }\rho ^{n-2k}a_{n,2k}H_{n-2k}\left( y|q\right) .
\end{equation*}
Let $S_n$ be the Hankel matrix of moments $m_k(y)$,
\begin{equation*}
S_{n}(y|q,\rho)= \left[ {\begin{array}{*{20}c} {m_0(y)} & {m_1(y)} & {\hdots} &
m_{n - 1}(y) \\ {m_1(y)} & {m_2(y)} & & {} \\ {\vdots} & {} & {\ddots} &
{\vdots} \\ m_{n - 1}(y) & {} & {\hdots} & m_{2n - 2}(y) \\ \end{array}} %
\right].
\end{equation*}
It is well known that $\det S_n$ is the product of the coefficients at the
third term of (\ref{+}), which implies the following.

\begin{corollary}
\label{C4} \label{wyznaczniki}$\det S_{n+1}/\det
S_{n}=[n]_q!\prod_{i=1}^{n}\left( 1-\rho^2 q^{i-1}\right) .$
\end{corollary}

Our second Hankel matrix has an even simpler form. As indicated in \cite{Ismail97}, 
\cite{IS03} the $q$-Hermite polynomials can be viewed as moments of a signed
measure, $H_n(x|q)=\int u^n \mu(du|x,q)$.
It turns out that 
measure $\mu(du|x,q)$ cannot be
positive even for a single value of $x$. 
To see this, consider the following $n\times n$ matrices
\begin{equation*}
M_n(x|q)=\left[ {\begin{array}{*{20}c} {H_0(x|q)} & {H_1(x|q)} & {H_2(x|q)}
& {\hdots} & {H_{n - 1}(x|q)} \\ {H_1(x|q)} & {H_2(x|q)} & {} & {} &
{H_{n}(x|q)} \\ {H_2(x|q)} & {} & {\ddots} & {} & {{H_{n+1}(x|q)}} \\
{\vdots} & {} & {} & {\ddots} & {\vdots} \\ {H_{n - 1}(x|q)} &
H_{n }(x|q) & {} & {\hdots} & H_{2n-2}(x|q) \\ \end{array}}
\right].
\end{equation*}

The following $q$-generalization of \cite[(3.55)]{Krattenthaler99} shows
that the determinants $\det M_n(x|q)$ are free of the variable $x$, and take negative values.

\begin{theorem}
\label{wyzn-Hermitow}
\begin{equation*}
\frac{\det M_{n+1} }{\det M_n}=(-1)^{n} q^{n(n-1)/2}[n]_q!.
\end{equation*}
\end{theorem}

\begin{proof}
Using (\ref{H}), we row-reduce the first column of the matrix. Namely, from
the second row of $M_{n+1}$, we subtract the first one multiplied by $x$.
Similarly, for $i\geq 3$, we subtract $x$ times row $i-1$ and add the $i-2$%
-th row multiplied by $[i-1]_{q}$. Taking (\ref{H}) into account, $\det
M_{n+1}(x|q)$ becomes
\begin{equation*}
\det \left[ {\begin{array}{*{20}c} {H_0} & {H_1} & {H_2} & {\hdots} &
{H_{n}} \\ 0 & {([0]-[1])H_0} & {([0]-[2])H_1} & {} & {([0]-[n])H_{n-1}} \\
0 & {([1]-[2])H_1} & {([1]-[3])H_2} & {} & {([1]-[n+1])H_{n}} \\
\vdots & {} & {} & {\ddots} & {\vdots} \\ 0 & ([n-1]-[n])H_{n-1} &
([n-1]-[n+1])H_{n} & {\hdots} & ([n-1]-[2n-1])H_{2n} \\
\end{array}}\right].
\end{equation*}%
Now, we use the fact that for $m\leq n$ we have
$
\lbrack n]_{q}-\left[ m\right] _{q}=q^{m}\left[ n-m\right] _{q}.
$
Thus $\det M_{n+1}(x|q)$ becomes
\begin{equation*}
\det \left[ {\begin{array}{*{20}c} {H_0} & {H_1} & {H_2} & {\hdots} & {H_{n
- 1}} \\ 0 & {-[1]H_0} & {-[2]H_1} & {} & {-[n]H_{n-1}} \\ 0 & {-q[1]H_1} &
{-q[2]H_2} & {} & {-q[n]H_{n}} \\ \vdots & {} & {} & {\ddots} & {\vdots} \\
0 & -q^{n-1}[1]H_{n-1} & -q^{n-1}[2]H_{n} & {\hdots} & -q^{n-1}[n]H_{2n} \\
\end{array}}\right].
\end{equation*}%
Expanding $\det M_{n+1}$ with respect to the first column, and
factoring out the common factors $-q^{i-1}$ from the $i$-th row
and $[j]_{q}$ from the $j$-th column of the resulting matrix, we
get
$$\det M_{n+1} =\left( -1\right)
^{n}q^{\sum_{i=1}^{n-1}i}\prod_{j=1}^{n}[j]\det M_{n} =\left(
-1\right) ^{n}q^{n\left( n-1\right) /2}[n]_{q}!\det M_{n}.$$
\end{proof}


Formula stated in Corollary \ref{C4} was originally discovered through
symbolic computations and motivated this paper. We were unable to find a
direct algebraic proof along the lines of the proof of Theorem \ref%
{wyzn-Hermitow}, and our search for the explanation of why $\det S_{n}(y)$
does not depend on $y$ lead us to Al-Salam--Chihara polynomials and identity
(\ref{MI}).

The the fact that Hankel determinants formed of certain linear combinations
of the $q$-Hermite polynomials do not depend on the argument of these
polynomials as exposed in Theorem \ref{wyzn-Hermitow} and Corollary \ref{C4}
is striking and unexpected to us. A natural question arises whether other
linear combinations have this property.

\medskip
{\bf Acknowledgement}
Part of the research of WB was conducted while visiting
Warsaw University of Technology. The authors thank
Jacek Weso{\l}owski for several helpful discussions.


\end{document}